\documentstyle[11pt]{article}
\newlength{\standardunitlength}
 \newtheorem{lemma}{Lemma}
\newtheorem{theorem}{Theorem} \newtheorem{prop}{Proposition}
\newenvironment{proof}{\noindent {\sc Proof:}}{$\Box$ \vspace{2 ex}}

\begin{document}

\begin{center}
Descent Algebras, Hyperplane Arrangements, and Shuffling Cards
\end{center}

\begin{center}
By Jason Fulman
\end{center}

\begin{center}
Dartmouth College
\end{center}

\begin{center}
Department of Mathematics
\end{center}

\begin{center}
Fulman@Dartmouth.Edu
\end{center}

\begin{abstract}
	Two notions of riffle shuffling on finite Coxeter groups are
given: one using Solomon's descent algebra and another using random
walk on chambers of hyperplane arrangements. These
coincide for types $A$,$B$,$C$, $H_3$, and rank two groups. Both
notions have the same, simple eigenvalues. The hyperplane definition
is especially natural and satisfies a positivity property when $W$ is
crystallographic and the relevant parameter is a good prime. The
hyperplane viewpoint suggests interesting connections with Lie theory
and leads to a notion of riffle shuffling for arbitrary real
hyperplane arrangements and oriented matroids. Connections with
Cellini's descent algebra are given. \end{abstract}

\begin{center}
\end{center}

\newpage

\section{Introduction and Background}

		Using ideas from \cite{BaD}, \cite{BB}, and
particularly \cite{BBHT}, we give a definition of card-shuffling
measures $M_{W,x}$ on finite Coxeter groups.  Here $x \neq 0$ is a
real number, and $M_{W,x}$ satisfies the measure property $\sum_{w \in
W} M_{W,x}(w) = 1$. In general these measures will be signed, i.e.  it
is possible that $M_{W,x}(w)<0$ for some element $w$. The measures
$M_{W,x}$ convolve and have nice eigenvalues.

		The measures $M_{W,x}$ have received considerable
attention in the cases that $W$ is of type $A$ or $B$. The type $A$
case appeared in \cite{BaD} in the theory of riffle shuffling. Results
for biased shuffles appear in \cite{F3}. Type $B$ riffle shuffles are
considered in \cite{BB}. As is evident from \cite{BB} and \cite{Ha},
$M_{A_n,x}$ and $M_{B_n,x}$ are related to the
Poincar\'e-Birkhoff-Witt theorem and to splittings of Hochschild
homology. Section 3.8 of \cite{SS} describes the measures $M_{A_n,x}$
in the language of Hopf algebras.

	We then give a second definition $H_{W,x}$ of riffle shuffling
for finite Coxeter groups using work of Bidigare, Hanlon, and Rockmore
\cite{BHR}. They describe how putting non-negative weights summing to
one on the faces of a hyperplane arrangement induces a random walk on
its chambers. This procedure will be recalled in Section
\ref{hyp1}. Brown and Diaconis \cite{BrD} generalize these chamber
walks to oriented matroids and give many examples. Bidigare \cite{B}
realized type $A$ shuffling as a special case, with face weights as
binomial coefficients. We give a definition for all finite Coxeter
groups using group theoretic weights. As a by-product, it will be seen
that if $W$ is crystallographic and $p$ is a good prime, then
$M_{W,p}(w) \geq 0$ for all $w \in W$. It will emerge that the
eigenvalues and multiplicites are the same as for the descent algebra
definition. The weights in the construction of $H_{W,x}$ are then
expressed in a completely combinatorial way, yielding a definition for
arbitrary real hyperplane arrangements.

	The follow-up work \cite{F1},\cite{F2} connects the measures
$H_{W,x}$ with the semisimple orbits of the adjoint action of a finite
group of Lie type on its Lie algebra. The paper \cite{F4} develops
analogous yet combinatorially quite different ideas for the case of
semisimple conjugacy classes of finite groups of Lie type, using a
probability measure $x_k$ on $W$ arising from Cellini's descent
algebra \cite{Ce1},\cite{Ce2}. Section \ref{Cel} of this paper ties
this measure in with $M_{W,k}$ and $H_{W,k}$ in type $C$ for odd
characteristic, and gives relations with physical models of card
shuffling in types $A$ and $C$.

	One long term goal, partially realized in \cite{F1},\cite{F4}
is to naturally associate to a semisimple conjugacy class or adjoint
orbit an element $w$ of the Weyl group, refining the well-known map
$\Phi$ to conjugacy classes of $W$. Letting $q$ be the size of the
base field, choosing a class (resp. orbit) uniformly at random and
applying the refined map gives $x_q$ (resp. $H_{W,q}$), at least in
some types. To say why refining $\Phi$ may be important, consider a
simple degree $n$ algebraic extension of $Q$ with minimal polynomial
$f(x)$ over $Z$. Reducing $f(x)$ mod a prime $p$ gives a semisimple
conjugacy class $c$ of $GL(n,p)$. At unramified primes, the Frobenius
automorphism, viewed as a permutation of the roots of $f$, has cycle
structure equal to that of $\Phi(c)$ (i.e. the number of $i$ cycles is
the number of degree $i$ irreducible factors of $f$ mod $p$). The
hope, and further motivation for this work, is that refining $\Phi$
will give refined number theoretic constructions.

\section{Descent Algebra Definition of Shuffling for Coxeter Groups}
\label{de}

	Let $W$ be a finite Coxeter group with $\Pi$ a base of
fundamental roots. For $w \in W$, let $Des(w)$ be the set of simple
positive roots mapped to negative roots by $w$ (also called the
descent set of $w$). For $J \subseteq \Pi$, let $X_J=\{w \in W| Des(w)
\cap J = \emptyset \}$ and $x_J = \sum_{w \in X_J} w$. Let $\lambda$
be an equivalence class of subsets of $\Pi$ under the equivalence
relation $J \sim K$ whenever $w(J)=K$ for some $w \in W$. Let
$\lambda(K)$ denote the equivalence class of $K$. For $K \subseteq J
\subseteq \Pi$ define $\mu_K^J = \frac{|\{w \in X_J: w(K) \subseteq
\Pi \}|}{|\lambda(K)|}$. Set $\mu_K^J=0$ if $K \not \subseteq J$. Let
$(\beta_K^J)$ be the matrix inverse of $(\mu_K^J)$. Define
$e_{\lambda}$ in the descent algebra of $W$ by

\[ e_{\lambda}  =  \sum_{J \in \lambda} \frac{\sum_{K \subseteq J} \beta_K^J
x_K}{|\lambda|}. \] These $e_{\lambda}$ are orthogonal idempotents which sum to
the identity element of $W$ \cite{BBHT}. Let $\| \lambda \|$ denote $|J|$ for any
$J$ in the equivalence class $\lambda$. 

{\bf Definition 1:} For $W$ a finite Coxeter group and $x \neq 0$,
define a signed probability measure $M_{W,x}$ on $W$ by

\[ M_{W,x} = \sum_{\lambda} \frac{e_{\lambda}}{x^{\| \lambda \|}}. \] For $w \in
W$, let $M_{W,x}(w)$ be the coefficient of $w$ in $M_{W,x}$.

\vspace{.5mm}
 
	Since the $e_{\lambda}$ are orthogonal idempotents, $M_{W,x}$ satisfies the
following convolution property in the group algebra of $W$:

\[ M_{W,x}M_{W,y} = M_{W,xy}. \]
\begin{prop} \label{signedmeasure} $M_{W,x}$ is a signed probability measure on
$W$.
\end{prop}

\begin{proof}
	Writing each $e_{\lambda}$ as $\sum_{w \in W} c_{\lambda}(w) w$ it must be
proved that 

\[ \sum_{w,\lambda} \frac{c_{\lambda}(w)}{x^{\| \lambda \|}} = 1. \]
This follows from the stronger assertion that $\sum_{w}
c_{\lambda}(w)$ is $0$ if $\| \lambda \|>0$ and is $1$ if $\| \lambda
\|=0$. Since $\beta_{\emptyset}^{\emptyset} = \frac{1}{|W|}$, it
follows by definition that $e_{\emptyset} = \frac{\sum_{w \in W}
w}{|W|}$. Thus $\sum_w c_{\lambda}(w) = 1$ if $\| \lambda \|=0$. Since
the $e_{\lambda}$ are idempotents, the value of $\sum_{w}
c_{\lambda}(w)$ is either 0 or 1. Since $\sum_{\lambda}
e_{\lambda}=1$, clearly $\sum_{w,\lambda} c_{\lambda}(w)=1$. Combining
this with the fact that $\sum_w c_{\lambda}(w) = 1$ if $\| \lambda
\|=0$ shows that $\sum_w c_{\lambda}(w) = 0$ if $\| \lambda \|>0$.
\end{proof}

	Proposition \ref{G2} illustrates the definition by computing
$M_{G_2,x}$.

\begin{prop} \label{G2}

\[ M_{G_2,x}(w) =             \left\{ \begin{array}{ll}
																																				\frac{(x+5)(x+1)}{12x^2} & \mbox{if
$d(w)=0$}\\
																																				\frac{(x+1)(x-1)}{12x^2} & \mbox{if
$d(w)=1$}\\
																																				\frac{(x-1)(x-5)}{12x^2}		& \mbox{if
$d(w)=2$}
																																				\end{array}
			\right.			 \]

\end{prop}

\begin{proof}
	Letting $V$ be the hyperplane in $R^3$ consisting of vectors whose coordinates
add to 0, it is well known that a root system consists of $\pm (\varepsilon_i -
\varepsilon_j)$ for $i<j$ and $\pm (2 \varepsilon_i - \varepsilon_j -
\varepsilon_k)$ where $\{i,j,k\}=\{1,2,3\}$. Let $A = \varepsilon_1 -
\varepsilon_2$ and $B=-2 \varepsilon_1+\varepsilon_2+\varepsilon_3$ be a base of
positive simple roots. All equivalence classes $\lambda$ of subsets of $\Pi$ have
size one. Some computation gives that

\begin{eqnarray*}
e_{\emptyset} & = & \frac{1}{12} x_{\emptyset}\\
e_A & = & -\frac{1}{4} x_{\emptyset} + \frac{1}{2}x_A\\
e_B & = & -\frac{1}{4} x_{\emptyset} + \frac{1}{2}x_B\\
e_{A,B} & = & \frac{5}{12} x_{\emptyset} -\frac{1}{2}x_A - \frac{1}{2}x_B+ x_{A,B}
\end{eqnarray*} from which the result follows.
\end{proof}

	Left multiplication in the group algebra of $W$ by $M_{W,x}$
can be thought of as performing a signed random walk on $W$. The
corresponding transition matrix is a $|W|$ by $|W|$
matrix. Proposition \ref{eig} determines its eigenvalues and evaluates
$M_{W,x}(w)$ when $w$ is the identity or longest element. It was known
previously for types $A$ and $B$ (part 1 in \cite{Ha}, the rest in
\cite{BaD} and \cite{BB}).

\begin{prop} \label{eig} Let $W$ be a finite Coxeter group of rank $n$. Let $id$
be the identity element of $W$ and $w_0$ the longest element of $W$. Let
$m_1,\cdots,m_n$ be the exponents of $W$.
\begin{enumerate}
\item The eigenvalues of $M_{W,x}$ are $\frac{1}{x^i}$ for $0 \leq i \leq n-1$
with multiplicity the number of $w \in W$ with fixed space of dimension $n-i$.
\item $M_{W,x}(id)=\frac{\prod_{i=1}^n (x+m_i)}{x^n|W|}.$
\item $M_{W,x}(w_0)=\frac{\prod_{i=1}^n (x-m_i)}{x^n|W|}.$
\end{enumerate}
\end{prop}

\begin{proof} The first claim is easily reduced to Theorem 7.15 of \cite{BBHT}.
For the second claim note from the first that summing over $\lambda$ with
$||\lambda||=i$ the coefficient of the identity in $e_{\lambda}$ gives
$\frac{1}{|W|}$ times the number of $w$ with fixed space of dimension $n-i$. Thus

\[ M_{W,x}(id)=\frac{1}{x^n |W|} \sum_{w \in W} x^{dim(fix(w))} =
\frac{\prod_{i=1}^n (x+m_i)}{x^n|W|} ,\] the second equality being an identity of
Shephard and Todd \cite{ST}. The third claim follows analogously to the second, but in the
proof of Theorem 7.15 of \cite{BBHT} one takes inner products with the alternating
representation instead of with the trivial representation.
\end{proof}

\section{Hyperplane Walk Definition of Shuffling for Coxeter Groups}
\label{hyp1}

	First it is necessary to review the paper \cite{BHR}. Let
$\cal{A}$ $= \{H_i : i \in I\}$ be a central hyperplane arrangement
(i.e. $\cap_{i \in I} H_i = 0$) for a real vector space $V$. Let
$\gamma$ be a vector in the complement of $\cal{A}$. Every $H_i$
partitions $V$ into three pieces: $H_i^0=H_i$, the open half space
$H_i^+$ of $V$ containing $\gamma$, and the open half space $H_i^-$ of
$V$ not containing $\gamma$. The faces of $\cal{A}$ are defined as the
non-empty intersections of the form $\cap_{i \in I} H_i^{\epsilon_i}$,
where $\epsilon_i \in \{0,-,+\}$. Equivalently, $\cal{A}$ cuts $V$
into regions called chambers and the faces are the faces of these
chambers viewed as polyhedra. A random process on chambers (henceforth
called the BHR walk) is then defined as follows. Assign weights $v(F)$
to the faces of $\cal{A}$ in such a way that $v(F) \geq 0$ for all $F$
and $\sum_F v(F)=1$. Pick a starting chamber $C_0$. At step $i$, pick
a face $F_i$ with the chance that $F_i=F$ equal to $v(F)$ and define $C_i$
to be the unique chamber whose closure contains $F_i$ and which among
such chambers is separated from $C_{i-1}$ by the fewest number of hyperplanes.

	For the remainder of this section, $\cal{A}$ will be the
arrangement of root hyperplanes for a finite Coxeter group $W$.

\begin{prop} (\cite{H}) The chambers of $\cal{A}$ correspond to the
elements of $W$. The faces of
$\cal{A}$ correspond to left cosets of parabolic subgroups of $W$. The
faces contained in the closure of $w$ are the left cosets $wW_J$.
\end{prop}

	The next lemma will be of use.

\begin{lemma} \label{descent} Let $C_0$ be the chamber of $\cal{A}$ corresponding
to the identity. Let $v(F)$ be the weight on a face $F$ in a BHR walk on a
finite Coxeter group. Then the chance that the chamber $C_1$ corrseponds to $w$
is equal to

\[ \sum_{K \subset \Pi-Des(w)} v(wW_K). \]

\end{lemma}

\begin{proof}
	The chance that $C_1$ corresponds to $w$ is equal to $\sum_F
v(F)$, where the sum is over all faces $F$ which are adjacent to $w$
and such that $w$ is the chamber adjacent to $F$ which is closest to
the identity. The faces adjacent to $w$ are the cosets $wW_K$, for $K$
an arbitrary subset of $\Pi$. The chambers adjacent to the face $wW_K$
are the elements of the coset $wW_K$. Proposition 1.10 of \cite{H}
shows that $w$ is the unique shortest element in the coset $wW_K$
precisely when $K \subset \Pi-Des(w)$.  \end{proof}

	Let $L$ be the set of intersections of the hyperplanes in
$\cal{A}$, taking $V \in L$. (This lattice is not the same as the face
lattice). Partially order $L$ by reverse inclusion. Recall that the
Moebius function $\mu$ is defined by $\mu(X,X)=1$ and $\sum_{X \leq Z
\leq Y} \mu(Z,Y)=0$ if $X<Y$ and $\mu(X,Y)=0$ otherwise. The
characteristic polynomial of $L$ is defined as

	\[ \chi(L,x) = \sum_{X \in L} \mu(V,X) x^{dim(X)}. \] (This is
not the standard definition in which the exponent of $x$ is
$n-rank(X)$. This distinction is important as the definitions do not
agree for subposets). For $J \subseteq \Pi$, let $Fix(W_J)$ denote the
fixed space of the parabolic subgroup $W_J$ in its action on $V$. Let
$L^{Fix(W_J)}$ be the subposet $\{Y \in L$$ |Y \geq
Fix(W_J)\}$. Finally, $N_{G_1}(G_2)$ denotes the normalizer of the
group $G_2$ in $G_1$.

{\bf Definition 2:} Let $W$ be a finite Coxeter group of rank $n$. Define 
$H_{W,x}(w)$ by the formula
\[ H_{W,x}(w) =
\sum_{K
\subseteq
\Pi-Des(w)}
\frac{|W_K|
\chi(L^{Fix(W_K)},x)}{x^n |N_W(W_K)||\lambda(K)|}.\]

\vspace{.5mm}

	Let $W$ be a finite irreducible crystallographic Coxeter group. Recall that $p$
is said to be a bad prime if $p$ divides the coefficient of some root of $W$ when
expressed as a linear combination of simple roots. Alternatively, $p$ is a bad
prime if it is less than the maximum exponent of $W$ but not equal to an exponent
of $W$. A prime $p$ is said to be good if $p$ is not bad.

\begin{prop} \label{pos} Let $W$ be a finite irreducible crystallographic Coxeter
group of rank $n$. Then if $p$ is a good prime for $W$, the measure $H_{W,p}$ can
be viewed as a special instance of the BHR walks with face weights
\[ v_p(wW_K) = \frac{|W_K| \chi(L^{Fix(W_K)},p)}{p^n |N_W(W_K)||\lambda(K)|}.\]
In particular, $H_{W,p}(w) \geq 0$ for all $w \in W$.
\end{prop}

\begin{proof} It will be proved in Section \ref{generalization} that
the face weights sum to 1. By Lemma \ref{descent}, it is sufficient to
show that $\chi(L^{Fix(W_K)},p)$ is non-negative for $p$ a good
prime. One of the main results of Orlik and Solomon \cite{OS} is the
factorization \[ \chi(L^{Fix(W_K)},x) = \prod_{i=1}^{dim(Fix(W_K))} (x
- b_i^K), \] where the $b_i^K$ are positive integers. From the tables
in their paper, the $b_i^K$ are all less than or equal to the maximum
exponent of $W$.  \end{proof}

{\bf Remarks:} \begin{enumerate} \item The BHR process can still be
considered with negative face weights, viewed as a transition matrix
with possibly negative entries. Using the face weights of Proposition
\ref{pos} where $W$ is crystallographic and $x$ is a bad prime gives a
first natural collection of examples where some of the face weights
are negative. As will emerge from Section \ref{generalization}, the
eigenvalues of the transition matrix are positive. We also note that
the $W$-invariance of the weights implies that the process can be
viewed as left multiplication by an element of the group algebra. It
would be interesting to prove that a cut-off phenomenon occurs for
these chains in the sense of \cite{BaD}.

\item Although the expression $H_{W,x}(w)$, viewed as a function of
$x$ for fixed $w$, factors into linear terms when $W$ is of type $A,B$
or $I$, this property does not hold in general (one simple
counterexample is taking $W$ of type $H_4$).  However the proof of
Corollary \ref{pos} implies that the face weights $v_{x}(wW_K)$ do
factor into linear terms as a function of $x$. This shows the
naturality of the hyperplane viewpoint. We observe (though it is not
clear from the definition) that in type $A$ these weights agree with
the binomial coefficients weights in \cite{BHR}.

\item The values $|\lambda(K)|$ have been tabulated \cite{C2}.
\end{enumerate}

\section{Riffle Shuffling for Real Hyperplane Arrangements}
\label{generalization}

	This section defines riffle shuffling for arbitrary real
hyperplane arrangements and explores some of its properties. As a
corollary it will be shown that $M_{W,x}$ and $H_{W,x}$ agree for many
types. In fact the main definition of this section extends to oriented
matroids by replacing the terms ``chamber'' and ``hyperplane
intersection'' by the terms ``tope'' and ``flat''. For clarity of
exposition we suppose that $\cal{A}$ is a real hyperplane arrangement,
possibly non-central, with a finite number of hyperplanes. Brown and
Diaconis \cite{BrD} verified that the BHR walks extend to this
setting.

	Recall the zero map from the face lattice to the intersection
lattice. This map $z$ sends a face $\cap_{i \in I} H_i^{\epsilon_i},
\epsilon_i \in \{0,-,+\}$ to $\cap_{i \in I: \epsilon_i=0}
H_i$. Geometrically, $z$ maps a face $F$ to its support, namely the
intersection of $V$ with all hyperplanes containing $F$.

\begin{lemma} \label{Zieg} (\cite{Z}) For all $Y \in L$,

\[ |z^{-1}(Y)| = \sum_{Z \geq Y} |\mu(Y,Z)| = \sum_{Z \geq Y} (-1)^{dim(Y)-dim(Z)}
\mu(Y,Z)=(-1)^{dim(Y)} \chi(L^Y,-1). \]

\end{lemma}

\begin{lemma} \label{O/S} (\cite{OS}) For all $X \in L$,

\[ \sum_{Y \in L \atop Y \geq X} \chi(L^Y,x) = x^{dim(X)}. \]

\end{lemma}
\begin{proof}
\begin{eqnarray*}
\sum_{Y \in L \atop Y \geq X} \chi(L^Y,x) & = & \sum_{Y \geq X} \sum_{Z \geq Y} \mu(Y,Z) x^{dim(Z)}\\
& = & \sum_{Z \in L} \sum_{Y: X \leq Y \leq Z} \mu(Y,Z) x^{dim(Z)}\\
& = & x^{dim(X)}.
\end{eqnarray*}
\end{proof}

{\bf Definition 3:} Define a one-parameter family of card-shuffling walks
on the chambers of a real hyperplane arrangement $\cal{A}$ in $n$ dimensions as
the BHR walks with face weights

\[ v_x(F) =  \frac{\chi(L^{z(F)},x)}{x^n |z^{-1}(z(F))|} = 
(-1)^{dim(z(F))} \frac{\chi(L^{z(F)},x)}{x^n \chi(L^{z(F)},-1)}. \]

	Lemmas \ref{Zieg} and \ref{O/S} imply that these face weights sum to one:

\[ \sum_F v_x(F) = \sum_{Y \in L \atop Y \geq V} \sum_{F:z(F)=Y}
\frac{\chi(L^{z(F)},x)}{x^n |z^{-1}(z(F))|} = \sum_{Y \in L \atop Y
\geq V} \frac{\chi(L^Y,x)}{x^n} =1.\]

 The $v_x(F)$ are not necessarily positive, in which case the BHR
walks take the extended meaning in the remark after Proposition
\ref{pos}. To show that Definition 3 extends Definition 2, the
following lemma will be helpful.

\begin{lemma} \label{O/S2} (\cite{OS})

\[ \frac{|N_W(W_K)||\lambda(K)|}{|W_K|} = (-1)^{dim(z(wW_K))}
\chi(L^{z(wW_K)},-1). \]

\end{lemma}

We remark that combining the fact that the face weights sum to one
with Lemma
\ref{O/S2} yields the possibly new identity

\[ \sum_{K \subseteq \Pi} (-1)^{n-|K|} \frac{|W|}{|W_K|}
\frac{\chi(L^{Fix(W_K)},x)}{\chi(L^{Fix(W_K)},-1)} = x^n. \] Setting $x=-1$ in
this identity gives the alternating sum formula

\[ \sum_{K \subseteq \Pi} (-1)^{|K|} \frac{|W|}{|W_K|} = 1, \] which
has a topological proof \cite{So1} as well as applications in the
invariant theory of Coxeter groups \cite{H}. Generalizations and
$q$-analogs related to the Steinberg character appear in Section 6.2
of \cite{C2}.

\begin{theorem} \label{specialize} The measures $H_{W,x}$ arise from Definition 3
with $\cal{A}$ equal to the arrangement of root hyperplanes of $W$.
\end{theorem}

\begin{proof}
	The theorem follows by showing that the face weights of Definition 2
are equal to the weights Definition 3 associates to the arrangement of root
hyperplanes of $W$. Thus it is necessary to prove that

\[ \frac{|W_K| \chi(L^{Fix(W_K)},x)}{x^n |N_W(W_K)||\lambda(K)|} =
(-1)^{dim(z(wW_K))} \frac{\chi(L^{z(wW_K)},x)}{x^n
\chi(L^{z(wW_K)},-1)}. \] Since $z(W_K)=Fix(W_K)$, $L^{z(wW_K)}$ is
isomorphic to $L^{Fix(W_K)}$, and thus the result follows from Lemma
\ref{O/S2}.  \end{proof}

	Theorem \ref{specialize} shows that the measures $H_{W,x}$ can
be easily computed from the tables of Orlik and Solomon. For example
one can check that $H_{W,x}$ agrees with the formula for $M_{W,x}$ on
the identity and longest element of $W$ (see the follow-up \cite{F1}
for a more conceptual proof). The point of including Definition 2 was
to show that the face weights can be expressed group theoretically.

\begin{theorem} The Markov chain associated to the hyperplane
arrangement $\cal{A}$ by Definition 1 has eigenvalues $\frac{1}{x^i}$
with multiplicity $\sum_{X \in L: dim(X)=n-i} |\mu(V,X)|$ for $0 \leq
i \leq n-1$. In particular, the eigenvalues and multiplicities of
$H_{W,x}$ agree with those of $M_{W,x}$.  \end{theorem}

\begin{proof}
	For the first assertion, \cite{BrD} (extending the paper
\cite{BHR} to the affine case) proves that the eigenvalues for their chamber walks
are indexed by elements
$X
\in L$ and are equal to
$\sum_{F:z(F)
\geq X} v(F)$ with multiplicity $|\mu(V,X)|$. Definition 3 and Lemma
\ref{O/S} imply that 

\[ \sum_{F:z(F) \geq X} v_x(F) =  \frac{1}{x^n} \sum_{Y \geq X} \chi(L^Y,x) = 
\frac{1}{x^{n-dim(X)}} .\] The second assertion follows from the factorization of the
characteristic polynomial of a root arrangement $\chi(L,x)=\prod_i (x-m_i)$
together with the result of Shephard and Todd used in part 2 of Proposition
\ref{eig}.
\end{proof}

	Theorem \ref{agree} shows that in many cases $H_{W,x}$ and
$M_{W,x}$ agree.

\begin{theorem} \label{agree} Let $W$ be a finite irreducible Coxeter group of
type $A,B,C,H_3$ or rank 2. Then $H_{W,x}=M_{W,x}$.
\end{theorem}

\begin{proof} The rank 2 cases are straightforward. For the case of $H_3$,
letting $d(w)$ be the number of descents of $w$ it follows from computations in
\cite{B3} that

\[ M_{H_3,x}(w) =             \left\{ \begin{array}{ll}
																																				\frac{(x+9)(x+5)(x+1)}{120x^3} & \mbox{if
$d(w)=0$}\\
																																				\frac{(x+5)(x+1)(x-1)}{120x^3} & \mbox{if
$d(w)=1$}\\
						                              \frac{(x+1)(x-1)(x-5)}{120x^3} & \mbox{if
$d(w)=2$}\\
																																				\frac{(x-1)(x-5)(x-9)}{120x^3}		& \mbox{if
$d(w)=3$}
																																				\end{array}
			\right.			 \] This checks with $H_{H_3,x}$ (which as explained earlier is
directly computable from tables in \cite{OS}).

	For the symmetric group, from \cite{BaD} it emerges that
$M_{S_n,x}(w)=\frac{{x+n-1-d(w) \choose n}}{x^n}$. To compute
$H_{S_n,x}$, Proposition 2.1 of \cite{OS} shows that
$\chi(L^{z(wW_K)},x)$ is equal to $(x-1)\cdots(x-(n-|K|)+1)$. Thus the
face weight Definition 3 associates to $W_K$ is $\frac{{x \choose
n-|K|}}{x^n}$, agreeing with the face weight for type $A$ shuffling in
\cite{BHR}. For type $B$, Proposition 2.2 of \cite{OS} shows that
$\chi(L^{z(wW_K)},x)$ is equal to $(x-1)(x-3)\cdots
(x-2(n-|K|)+1)$. Thus

\begin{eqnarray*}
H_{B_n,x}(w)&=& \sum_{K \subseteq \Pi-Des(w)}
\frac{(-1)^{dim(z(W_K))} \chi(L^{z(wW_K)},x)}{x^n
\chi(L^{z(wW_K)},-1)}\\
& = & \frac{1}{x^n} \sum_{j=0}^{n-d(w)}
{n-d(w) \choose j} {\frac{x-1}{2} \choose n-j}\\
& = & \frac{{\frac{x-1}{2}+n-d(w) \choose n}}{x^n},
\end{eqnarray*} agreeing
with the formula for $M_{B_n,x}$ in \cite{BB}. The argument for type
$C$ is identical.  \end{proof}

	Theorem \ref{agree} implies that the measures $H_{W,x}$
convolve when $W$ is of type $A,B,C,H_3$ or a rank 2 group.

\section{Cellini's Descent Algebra} \label{Cel}

	To begin we define for any Weyl group $W$ and positive integer
$k$ an element $x_k$ of the group algebra of $W$. This construction is
due to Cellini \cite{Ce1} (the definition which follows differs
slightly from hers: it is inverse, it uses her Corollary 2.1, and it
is renormalized so as to yield a probability measure).

	Letting $\alpha_0$ denote the highest root, let
$\tilde{\Pi}=\Pi \cup \alpha_0$. Define the cyclic descent $Cdes(w)$
to be the elements of $\tilde{\Pi}$ mapped to negative roots by $w$,
and let $cd(w)=|Cdes(w)|$. Let $Y$ be the coroot lattice. Then define
$a_{k,I}$ by

\[ \left\{ \begin{array}{ll}
 |\{ t \in Y| <\alpha_0,t> = k, <\alpha_i,t> = 0 \ for \ \alpha_i
\in I-\alpha_0, <\alpha_i,t> >0 \ for \ \alpha_i \in \tilde{\Pi}-I\}| &
\mbox{if $\alpha_0 \in I$}\\
 |\{ t \in Y| <\alpha_0,t> < k, <\alpha_i,t> = 0 \ for \ \alpha_i \in
I, <\alpha_i,t> >0 \ for \ \alpha_i \in \Pi-I\}| &
\mbox{if $\alpha_0 \not \in I$} \end{array}\right.\] Finally, define an
element $x_k$ of the group algebra of $W$ by \[ x_k = \sum_{w \in W}
\left(\frac{1}{k^r} \sum_{I \subseteq \tilde{\Pi} - Cdes(w)} a_{k,I} \right) w.\] From Cellini (loc. cit.), it follows that the $x_k$ satisfy
the following two desirable properties:

\begin{enumerate} \item (Measure) The sum of the coefficients in the
expansion of $x_k$ in the basis of group elements is 1. In
probabilistic terms, the element $x_k$ defines a probability measure
on the group $W$.

\item (Convolution) $x_k x_h = x_{kh}$.
\end{enumerate}

	The above definition of $x_k$ is computationally convenient
for this paper. We note that Cellini (loc. cit.) constructed the $x_k$
in the following more conceptual way. Let $W_k$ be the index $k^r$
subgroup of the affine Weyl group generated by reflections in the
hyperplanes corresponding to $\{\alpha_1,\cdots,\alpha_r\}$ and also
the hyperplane $\{<x,\alpha_0>=k\}$. There are $k^r$ unique minimal
length coset representatives for $W_k$ in the affine Weyl group, and
$x_k$ is obtained by projecting them to the Weyl group.

	In type $A$ the elements $x_k$ are very different from the
measures $M_{W,k}$ and $H_{W,k}$ considered in Sections \ref{de} and
\ref{hyp1} (this is not too suprising as semsimple conjugacy classes
and semisimple adjoint orbits, though equinumerous, are quite
different combinatorially). For example the coefficients of $w$ in
$x_k$ are not polynomials in $k$, but depend on the value of $k$ mod
$n$. Nevertheless, as Proposition \ref{coincide} shows, all three
notions coincide in type $C$ when $k$ is odd.

\begin{prop} \label{coincide} Let $W$ be of type $C$ and $k$ be an odd
positive integer. Then $x_k=M_{W,k}=H_{W,k}$. \end{prop}

\begin{proof} From Theorem 1 of \cite{Ce2}, it follows that if $w$ has
$d$ descents and $k$ is odd, then the coefficient of $w$ in $x_k$ is
\begin{eqnarray*} \frac{1}{k^n} \sum_{l=d(w)}^n {\frac{k-1}{2} \choose
l} {n-d(w) \choose l-d(w)} & = & \frac{1}{k^n} \sum_{l=d(w)}^n
{\frac{k-1}{2} \choose l} {n-d(w) \choose n-l}\\ & = & \frac{1}{k^n}
\sum_{l=0}^n {\frac{k-1}{2} \choose l} {n-d(w) \choose n-l}\\ & = &
\frac{1}{k^n} {\frac{k-1}{2}+n-d(w) \choose n}.  \end{eqnarray*} This
visibly agrees with the formulas for $M_{W,k}$ and $H_{W,k}$ in
Theorem \ref{agree}.
\end{proof}

	For those interested in physical models of card-shuffling, we
remark that the inverse of the shuffles in Proposition \ref{coincide}
admit the following physical description (observed for $k=3$ in
\cite{BB}, and different from the hyperoctahedral shuffles in
\cite{BaD}).

	Step 1: Choose $2k+1$ numbers $j_1,\cdots,j_{2k+1}$
multinomially with the probability of getting $j_1,\cdots,j_{2k+1}$
equal to $\frac{{n \choose j_1,\cdots,j_{2k+1}}}{(2k+1)^n}$. Make
$2k+1$ stacks of cards of sizes $j_1,\cdots,j_{2k+1}$
respectively. Then flip over the even numbered stacks.

	Step 2: Drop cards from packets with probability proportional
to packet size at a given time. Thus if the packets have size
$A_1,\cdots,A_{2k+1}$ at a given time, drop from packet $i$ with
probability $\frac{A_i}{A_1+\cdots+A_{2k+1}}$.

	Proposition \ref{physical} shows that $x_2$ has a simple
physical description in type $A$.
	
\begin{prop} \label{physical} When $W$ is the symmetric group
$S_{2n}$, $x_2^{-1}$ has the following probabilistic interpretation:

	Step 1: Choose an even number between $1$ and $2n$ with the
probability of getting $2j$ equal to $\frac{{2n \choose
2j}}{2^{2n-1}}$. From the stack of $2n$ cards, form a second pile of
size $2j$ by removing the top $j$ cards of the stack, and then putting
the bottom $j$ cards of the first stack on top of them.

	Step 2: Now one has a stack of size $2n-2j$ and a stack of
size $2j$. Drop cards repeatedly according to the rule that if stacks
$1,2$ have sizes $A,B$ at some time, then the next card comes from
stack $1$ with probability $\frac{A}{A+B}$ and from stack 2 with
probability $\frac{B}{A+B}$. (Equivalently, one chooses uniformly at
random one of the ${2n \choose 2j}$ interleavings preserving the
relative orders of the cards in each stack).

	The description of $x_2^{-1}$ is the same for the symmetric
group $S_{2n+1}$, except that at the beginning of Step 1, the chance
of getting $2j$ is $\frac{{2n+1 \choose 2j}}{2^{2n}}$ and at the
beginning of Step 2, one has a stack of size $2n+1-2j$ and a stack of
size $2j$. \end{prop}

\begin{proof} Suppose that $W$ is $S_{2n}$, the argument for
$S_{2n+1}$ being similar. The coroot lattice is all vectors with
integer components and zero sum with respect to a basis
$e_1,\cdots,e_{2n}$, that $\alpha_i=e_i-e_{i+1}$ for $i=1,\cdots,2n-1$
and that $\alpha_0=e_1-e_{2n}$. The elements of the coroot lattice
contributing to some $a_{2,I}$ are:

\[ \begin{array}{ll} (0,0,\cdots,0,0) & I=\tilde{\Pi}-\alpha_0\\
(1,0,0,\cdots,0,0,-1) & I=\tilde{\Pi}-\{\alpha_1,\alpha_{2n-1}\}\\
(1,1,0,0,\cdots,0,0,-1,-1) &
I=\tilde{\Pi}-\{\alpha_2,\alpha_{2n-2}\}\\ \cdots & \cdots\\
(1,1,\cdots,1,0,0,-1,\cdots,-1,-1) &
I=\tilde{\Pi}-\{\alpha_{n-1},\alpha_{n+1}\}\\
(1,1,\cdots,1,1,-1,-1,\cdots,-1,-1) & I=\tilde{\Pi}-\alpha_n
\end{array} \] One observes that the permutations in the above card
shuffling description all contribute to $u_I$ where

\[ I = \left\{ \begin{array}{ll} \tilde{\Pi}-\alpha_0 & \mbox{if
$2j=0$}\\ \tilde{\Pi}-\{\alpha_k,\alpha_{2n-k}\} & \mbox{if $2j=2
min(k,2n-k)$}\\ \tilde{\Pi}-\alpha_n & \mbox{if $2j=2n$}\end{array}
			\right.  \] As the total number of such
permutations is $\sum_{j=0}^m {2n \choose 2j}=2^{2n-1}$ and $x_2$ is a
sum of $2^{2n-1}$ group elements, the proof is complete. \end{proof}

	We plan to later study physical models of the various shuffles
in this paper for other finite Coxeter groups, using Eriksson's
natural permutation representations.

\section{Acknowledgements} The author thanks Persi Diaconis for
keeping him well informed of his work, and a careful referee and
reviewer for corrections. Part of this work was supported by an NSF
Postdoctoral Fellowship.


\begin{thebibliography}{AAA}

\bibitem [BaD]{BaD} Bayer, D. and Diaconis, P., Trailing the dovetail
shuffle to its lair. {\it Ann. of Appl. Probab.} {\bf 2}, Number 2
(1992), 294-313.

\bibitem [BB]{BB} Bergeron, F. and Bergeron, N., Orthogonal idempotents in the
descent algebra of $B_n$ and applications. {\it J. Pure Appl.
Algebra} {\bf 79}, Number 2 (1992), 109-129.

\bibitem [BBHT]{BBHT} Bergeron, F., Bergeron, N., Howlett, R.B., and Taylor, D.E.,
A decomposition of the descent algebra of a finite Coxeter group. {\it J.
Algebraic Combin.} {\bf 1} (1992), 23-44.

\bibitem [B3]{B3} Bergeron, F., and Bergeron, N., Symbolic manipulation for the
study of the descent algebra of finite Coxeter groups. {\it J. Symbolic Comput.}
{\bf 14} (1992), 127-139.

\bibitem [B]{B} Bidigare, P., Hyperplane arrangement face algebras and their
associated Markov chains, {\it Ph.D. Thesis}, University of Michigan, 1997.

\bibitem [BHR]{BHR} Bidigare, P., Hanlon, P., and Rockmore, D., A combinatorial
description of the spectrum of the Tsetlin library and its generalization to
hyperplane arrangements. To appear in {\it Duke Math. J.}.

\bibitem [BrD]{BrD} Brown, K., and Diaconis, P., Random walk and
hyperplane arrangements. {\it Ann. of Probab.} {\bf 26} (1998),
1813-1854.

\bibitem [C1]{C} Carter, R., {\it Finite groups of Lie type}. John
Wiley and Sons, 1985.

\bibitem [C2]{C2} Carter, R., Conjugacy classes in the Weyl
group. {\it Composito Math.} {\bf 25} (1972), 1-59.

\bibitem [Ce1]{Ce1} Cellini, P., A general commutative descent algebra. {\it
J. Algebra} {\bf 175} (1995), 990-1014.

\bibitem [Ce2]{Ce2} Cellini, P., A general commutative descent algebra II.
The Case $C_n$. {\it J. Algebra} {\bf 175} (1995), 1015-1026.

\bibitem [F1]{F1} Fulman, J., Semisimple orbits of Lie algebras and
card shuffling measures on Coxeter groups. Technical report.

\bibitem [F2]{F2} Fulman, J., Counting semisimple orbits of finite Lie
algebras by genus, {\it J. Algebra} {\bf 217} (1999), 170-179.

\bibitem [F3]{F3} Fulman, J., The combinatorics of biased riffle shuffles.
{\it Ann. of Combin.} {\bf 2} (1998), 1-6.

\bibitem [F4]{F4} Fulman, J., Cellini's descent algebra and semisimple
conjugacy classes of finite groups of Lie type. Technical report.

\bibitem [Ha]{Ha} Hanlon, P., The action of $S_n$ on the components of the Hodge
decompositions of Hochschild homology, {\it Michigan Math. J.} {\bf 37}
105-124.

\bibitem [H]{H} Humphreys, J., {\it Reflection groups and Coxeter groups}.
Cambridge Studies in Advanced Mathematics {\bf 29}, Cambridge University Press,
Cambridge.

\bibitem [OS]{OS} Orlik, P., and Solomon, L., Coxeter arrangements. {\it
Proc. Symposia in Pure Math.} {\bf 40} (1983), Part 2, 269-291.

\bibitem [SS]{SS} Shnider, S., and Sternberg, S., {\it Quantum groups.}
Graduate Texts in Mathematical Physics, {\bf II}. International Press, 1993. 

\bibitem [ST]{ST} Shephard, G.C., and Todd, J.A., Finite unitary
reflection groups {\it Canadian J. Math.} {\bf 6} (1954), 274-304.

\bibitem [So1]{So1} Solomon, L., The orders of the finite Chevalley groups. {\it J.
Algebra} {\bf 3} (1966), 376-393.

\bibitem [Z]{Z} Zaslavsky, T., Facing up to arrangements: face-count formulas for
partitions of space by hyperplanes. {\it Mem. Amer. Math. Soc.} {\bf 1} (154).

\end{thebibliography}
\end{document}